\documentclass[10pt]{article}
\usepackage[english]{babel}
\usepackage{color}
\usepackage[dvipsnames]{xcolor}
\usepackage{amsmath}
\usepackage{physics}
\usepackage{xr-hyper}
\usepackage[hidelinks]{hyperref}
\usepackage{diagbox}
\usepackage{multirow}
\usepackage{amsfonts}
\usepackage{amssymb}
\usepackage{latexsym}
\usepackage{graphicx}
\usepackage{adjustbox}
\usepackage{float}
\usepackage{xr}
\usepackage{caption}
\usepackage{enumitem}
\usepackage{subfigure}
\usepackage{seqsplit}
\usepackage{appendix}
\usepackage[normalem]{ulem}
\usepackage[margin=0.75in]{geometry}
\usepackage{ulem}
\usepackage{cancel}
\usepackage{xcolor}

\usepackage{mathtools}
\usepackage[nodisplayskipstretch]{setspace} 

\newcommand{\pa}{\partial}

\newcommand{\ta}{\tau}

\newcommand{\non}{\nonumber}

\newtheorem{thm}{Theorem}[section]

\newtheorem{rem}[thm]{Remark}

\newtheorem{assumption}[thm]{Assumption}

\newcommand{\Int}{\int\limits}

\usepackage[nottoc]{tocbibind}

\begin{document}

	\title{Second order exponential splittings in the presence of unbounded and time-dependent operators: construction and convergence}
	\author{ Karolina Kropielnicka\footnote{Institute of Mathematics, Polish Academy of Sciences, Warsaw, Poland.}\qquad Juan Carlos del Valle\footnote{Faculty of Mathematics, Physics, and Informatics, University of Gda\'nsk, Gda\'nsk, Poland.}}
	\date{\today}
	
	\maketitle
	
	
	\begin{abstract}
		For linear differential equations of the form $u'(t)=[A + B(t)] u(t)$, $t\geq0$, with a possibly unbounded operator $A$, we construct and deduce error bounds for two families of second-order exponential splittings.  
			The role of quadratures when integrating the twice-iterated Duhamel's formula is reformulated: we show that their choice defines the structure of the splitting. Furthermore, the reformulation allows us to consider quadratures based on the Birkhoff interpolation to obtain not only pure-stages splittings but also those containing derivatives of $B(t)$ and  commutators of $A$ and $B(t)$. In this approach, the construction and error analysis of the splittings are carried out simultaneously. We discuss the accuracy of the members of the families. Numerical experiments are presented to complement the theoretical consideration.
		
	\end{abstract}

	\section{Introduction}
	{\color{black}
		Linear evolutionary differential equations of the form
		\begin{equation}
			u'(t)  =[A+B(t)]u(t),\qquad u(0)=u_0,\qquad t\in[0,T],
			\label{main}
		\end{equation}
		describe a wide variety of phenomena in the context of biology,  engineering,  physics, and many others. For instance, a remarkable case is when $A$ and $B(t)$ are skew-hermitian operators, then (\ref{main}) describes phenomena in the framework of  quantum mechanics.
		
		Numerical approximate solutions for (\ref{main}) can be constructed using various approaches like implicit and explicit Runge-Kutta methods,  
		exponential splittings or  exponential integrators, to mention a few remarkable; see, e.g. \cite{BlanesCasaMurua,HairerBook,hochbruck_ostermann_2010,McLachlan} and references therein. In particular,  exponential integrators are rooted in the celebrated Duhamel's formula, also known as variation of constants\footnote{We use both names indistinctibly. }. An extensive review of this kind of exponential integrators can be found in \cite{hochbruck_ostermann_2010}. 
		Exponential splittings are, on the other hand,  based on the decomposition of the vector field, thus, they are  essentially compositions of semigroups generated by the linear operators $A$, $B(t)$ at fixed $t$, and/or their commutators; see. e.g., \cite{BlanesCasaMurua,McLachlan}. Both approaches have their own assets depending on  $A$, $B(t)$ and the initial condition.\footnote{Exponential integrators can be tailored for rough initial data, while splittings methods may require higher regularity. On the other hand, assuming skew-Hermitian $A$ and $B(t)$, splitting methods preserve the norm of the solution, which is a desirable feature in the context of quantum mechanics. In contrast,  exponential integrators are usually not  preserving the mass.}
		
		
		In the case of a bounded operator $A$ and time independent $B$, a wide variety of pure-stages splittings has been  applied to describe several phenomena, see \cite{Schlick,Tong}. Taylor's theorem, Baker--Campbell--Hausdorff formula and associated order conditions are the standard tools to perform the error analysis of the splittings, see \cite{BlanesCasaMurua,McLachlan}. Extension of these splittings  to the case of time-dependent $B(t)$ can  be done using the notion of a superoperator, see  \cite{Chin}, or by translating  (\ref{main}) to and autonomous equation defining $t$ as a new coordinate, see  \cite{BlanesBook}. Another approach used to tackle splittings for the case of time-dependent $B(t)$ may be based, for example, on the symmetric Fer expansion \cite{Zanna} or  a decomposition applied to the truncated Magnus expansion \cite{Iserles}.
		
		In the case of an unbounded operator $A$ (and time-independent $B$), the derivation and analysis of splitting methods require other considerations because semigroups generated by unbounded operators do not undergo Taylor expansions, see \cite{Pazy}. Additionally, the convergence of Magnus, Dyson or Fer expansions is not guaranteed, with the exception of the truncated fourth-order Magnus expansion \cite{Thalhammer}.
		The convergence of Strang splitting in the presence of unbounded $A$ was proved in the seminal work  \cite{Lubich}. From the approximate integration of an iterated Duhamel's formula, the authors remarkably claimed: ``\textit{A basic observation is that the principal
			error terms} (of the splitting) \textit{are just quadrature errors}."  Using a similar approach, these ideas were revisited and extended to an arbitrary $p$-stages splitting in \cite{Mechthild}. The latter work also shows that the main error terms come from quadratures used for the  approximations of the integrals in the iterated Duhamel's formula. 
		Specifically, the role of quadratures is summarized by the author in the following sentence ``\textit{we associate the order conditions with quadrature order conditions for multiple integrals}."
		
		In the present work, assuming $A$ to be unbounded and  $B(t)$ to be bounded and time-dependent,  we reformulate the role of the quadratures used in the iterated Duhamel's formula, showing how their choice influences the error and structure of the splitting (e.g. the number of stages and/or the appearance of commutators in the exponents). Thus, quadratures are the building blocks behind the splittings. More precisely, we show how the exponentials containing $A$, $B(t)$,  its commutators and time-derivatives may appear in the splittings by considering quadratures based on Birkhoff interpolation. We refer to such quadratures as Birkhoff quadratures.  
		\\
		We focus on two representative  second-order families of splittings characterized by real one-parameter $\tau$. Based on these families, the reformulation of quadratures is established and illustrated.
		Firstly, we consider a  family with \textit{pure} stages only,
		\begin{equation}\label{family}
			\mathcal{F}(h,\tau) :=  e^{h\ta A}e^{\frac{h}{2} B(h(1-\ta))}e^{h(1-2\ta) A}e^{\frac{h}{2}B(h\ta) }e^{h\ta A} ,\qquad \tau\in[0,\tfrac{1}{2}],
		\end{equation}
		where $\tau$ is a  parameter and  $h$ is the time step.\footnote{In the case of time-independent $B$, family $\mathcal{F}(h,\tau)$ has been presented in \cite{McLachlan-II,McLachlan-I} but the error analysis was not carried out.} 
		Three-stages splittings can  be  obtained only for $\tau=0$ and $1/2$, leading to midpoint Strang splittings.  Additionally, ${\tau=1/4}$ leads to the composition of two midpoint Strang splittings with half the step. 
		Secondly, we consider the three-stages family 
		\begin{equation}
			\label{familyII}
			\mathcal{D}(h,\tau):= e^{h(1-\tau)A}e^{hB(h\tau)+\frac{h^2(1-2\tau)}{2}\mathcal{C}(h\tau)}e^{h\tau A}\ ,\quad\tau\in[0,1],
		\end{equation}
		where $\mathcal{C}(t)=[B(t),A]+B'(t)$.
		In this case, the values $\tau=0,1$ lead to two stages, and we revisit the midpoint Strang splitting when $\tau=1/2$.  In contrast to (\ref{family}), the family $\mathcal{D}(h,\tau)$ stands out due to the presence of the commutator $[B(t),A]$ and $B'(t)$. In fact, $\mathcal{D}(h,\tau)$ is the minimal example of a splitting  containing a commutator. The error analysis of this kind of splitting is out of the scope of \cite{Mechthild}.
		\\
		
		Despite the different structure of $\mathcal{F}(h,\tau)$ and $\mathcal{D}(h,\tau)$, we show how the same elementary arguments based on quadratures can be applied to their construction and error analysis. Thus, the present work can be regarded as the first step towards a generalization of the previous works \cite{Lubich} and \cite{Mechthild} embracing also  splittings featured by exponentials containing  time derivatives of $B(t)$ and commutators of $A$ and $B(t)$.
		\\
		
		This work has the following structure. In Section \ref{Family_F}, we present the derivation and error analysis of family $\mathcal{F}(h,\tau)$. More precisely, we formulate the basic assumptions on equation (\ref{main})  and derive  $\mathcal{F}(h,\tau)$ in Subsection \ref{derivation}, indicating all error terms that will be investigated in Subsection \ref{error}. The discussion on the role of $\tau$ is presented in Subsection \ref{accuracy}. In the same way, the family $\mathcal{D}(h,\tau)$ is derived and analyzed in Section \ref{Family_D} and its Subsections. Then, in Section \ref{example}, we present numerical examples based on the linear Sch\"rodinger and transport equations, where we investigate the influence of $\tau$ on the performance of  $\mathcal{F}(h,\tau)$ and $\mathcal{D}(h,\tau)$. Finally, our results are summarized in Section \ref{conclusions}.
	}
	\section{Family of Splittings $\mathcal{F}(h,\tau)$}\label{Family_F}
	
	\subsection{Derivation of the family  of splittings} \label{derivation}
	
	For clarity of exposition, we will derive the local error of the first step only, c.f. (\ref{family}), understanding that the whole procedure applies to any time step $t_n\in[0,T]$
	$$
	u(t_n+h)\ \approx \ e^{h\ta A}e^{\frac{h}{2} B((t_n+h)(1-\ta))}e^{h(1-2\ta) A}e^{\frac{h}{2}B((t_n+h)\ta) }e^{h\ta A}u(t_n) ,\qquad \tau\in[0,\tfrac{1}{2}],
	$$
	with minor and obvious modifications on the assumptions. 
	
	Throughout this paper, we consider problem (\ref{main}) in the abstract setting where $X$ is a Banach space with the norm and induced operator norm $\|\cdot\|$ and make the following assumption.
	
	\begin{assumption}\label{main_assum}
		We assume that densely defined and closed, linear (possibly unbounded) operator ${A:D(A)\subset X \rightarrow X}$ is a  generator of a strongly continuous $C_0$-semigroup ${\rm e}^{tA}$ on $X$, and that, for each fixed $s\in [0,h]$, $B(s)$ is a bounded, linear operator acting on $X$, that is $\forall_{s\in[0,h]}B(s)\in\mathfrak{B}(X)$, where $\mathfrak{B}(X)$ is the space of bounded linear operators defined on $X$. Moreover we assume that $B\in C^2([0,h],\mathfrak{B}(X))$ and seek solutions $u\in C^1([0,h],X)$.
	\end{assumption}
	
	For the convenience of notation let us define the space $Y:=C^1([0,h],X)$ with the standard norm \\ ${\|v(\cdot)\|_Y:=\max_{s\in[0,h]}\|v(s)\|}$. Given that both  $A$ and, for each fixed $s\in[0,h]$, $B(s)$ generate $C_0$-semigroups, there exists a constant $C_h$ depending on $h$ only, such that $\|{\rm e}^{tA}\|\leq C_h $ and $\|{\rm e}^{tB(\cdot)}\|_Y\leq C_h$  for $t\in[0,h]$.
	
	We wish to show that depending on the error constant, which may involve for example $[B(\cdot),A]$ or $[[B(\cdot),A],A]$, the (local) error of family of integrates (\ref{family}) scales like $\mathcal{O}(h^2)$ or $\mathcal{O}(h^3)$, respectively. For simplicity of the presentation, we derive the family of (\ref{family}) in two separate subsections according to the two desired order of error estimates. In both cases, the starting point is  expressing the solution of  (\ref{main}) via the well-known variation of constants formula,
	\begin{equation}
		u(h)= e^{hA}u_0\ +\ \int_0^h\,e^{(h-s)A}B(s)u(s) \,\text{d}s.
		\label{Duhamel}
	\end{equation}
	
	\subsubsection{First-order version}
	
	One iteration of (\ref{Duhamel}) leads to
	\begin{equation}\label{variation1}
		u(h)= e^{hA}u_0+ \int_0^he^{(h-t_1)A}B(t_1)e^{t_1A}u_0\,\text{d}t_1 + R_{V_1}, 
	\end{equation}
	with 
	\begin{equation}\label{RV1}
		R_{V_1}=\int_0^{h}\int_{0}^{t_1}e^{(h-t_1)A}B(t_1)e^{(t_1-t_2)A}B(t_2)u(t_2)\,\text{d}t_2\text{d}t_1.
	\end{equation}
	The next step is to approximate the integral of (\ref{variation1}). The choice of the quadratures is essential since  determines the structure of the splitting. For the family  $\mathcal{F}(h,\tau)$, we use
	\begin{equation}\label{quadrature1}
		I_1:=\  \frac{h}{2}f_1(h(1-\ta))\ +\ \frac{h}{2}f_1(h\ta),\qquad f_1(t_1)=e^{(h-t_1)A}B(t_1)e^{t_1A}u_0,
	\end{equation}	
	which for $\tau=\tfrac{1}{2}$  coincides with the midpoint rule, while for $\ta=0$ leads to the trapezoidal one. 
	Its  error 
	\begin{equation}\label{RI1}
		R_{I_1}:=\int_0^{h}f_1(t_1)\,\text{d}t_1\ -\ I_1, 
	\end{equation}	
	of order $\mathcal{O}(h^2)$ will be derived in Section \ref{error}.
	
	After incorporating the explicit form of $f_1(t_1)$ to $I_1$, see (\ref{quadrature1}),  the following equality holds
	\begin{align}
		e^{hA}u_0+I_1 = &\, {\rm e}^{h\tau A} \,  {\rm Id}\,  e^{h(1-2\tau) A} \,  {\rm Id} \,  {\rm e}^{h\tau A}u_0\non \\
		+&\frac{h}{2}\,e^{ h\tau A}\,  B(h(1-\tau))\,  e^{h(1-2\tau) A}\,  {\rm Id}\,  e^{h\tau A}u_0\non \\ 
		+&\frac{h}{2}\,e^{h\tau A}\,  {\rm Id}\,  e^{h(1-2\tau) A}\,  B(h\tau)\,  e^{h\tau A}u_0 \label{still_missing_1}.
	\end{align}
	Defining $R_{E_1}$ as 
	\begin{equation}
		R_{E_1}:=\,\frac{h^2}{4}e^{ h\tau A}\,  B(h(1-\tau))\,  e^{h(1-2\tau)A}\,  B(h\tau)\,  e^{h\tau A}u_0,\label{RE1}
	\end{equation} 
	we can observe that
	\begin{equation*}
		e^{hA}u_0+I_1+R_{E_1} =e^{ h\tau A}\left[{\rm Id}\, +\,  \frac{h}{2}\, B(h(1-\tau))\right]e^{h(1-2\tau)A}\left[{\rm Id}\, +\,  \frac{h}{2}\,B(h\tau)  \right]e^{h\tau A}u_0.
	\end{equation*}
	Given that the operator $B(s)$ is bounded at fixed $s\in[0,h]$, one can easily observe that expressions inside square parenthesis can be approximated  by the  semigroups  $e^{\frac{h}{2} B(h(1-\ta))}$ and $e^{\frac{h}{2} B(h\ta)}$, respectively. In this way, $R_{E_1}$ is related to the \textit{reconstruction} of exponentials (semigroups).  This approximation, however, introduces an additional error term
	\begin{align}
		R_{S_1}:=
		&
		-\frac{h^2}{8}\,e^{ h\tau A}\,  B^2(h(1-\tau))\,  e^{\frac{h}{2}\zeta_1B(h(1-\tau))}\,  e^{h(1-2\tau) A}\,  e^{\frac{h}{2} B(h\ta)}\,  e^{h\tau A}u_0\non \\
		&
		-\frac{h^2}{8}\,e^{h\tau A}\,  e^{\frac{h}{2} B(h(1-\ta))}\,  e^{h(1-2\tau) A}\,  B^2(h\tau)\,  e^{\frac{h}{2}\zeta_2B(h\tau)}\,  e^{h\tau A}u_0 \non \\
		&
		\,+\frac{h^4}{64}\,e^{h\tau A}\,  B^2(h(1-\tau))\,  e^{\frac{h}{2}\zeta_1B(h(1-\tau))}\,   e^{h(1-2\tau)A}\, B^2(h\tau)\,  e^{\frac{h}{2}\zeta_2B(h\tau)}\,  e^{h\tau A}u_0 \label{RS1},
	\end{align}
	for certain $\zeta_1\in[0,1]$ and $\zeta_2\in[0,1]$.
	As a result,  the connection between the exact solution of (\ref{main}) and the family of integrators (\ref{family}) reads 
	\begin{equation}\label{Error_1}
		u(h)\ = \ \mathcal{F}(h,\ta)u_0+R_{V_1}+R_{I_1}-R_{E_1}+R_{S_1},\qquad \tau\in[0,\tfrac{1}{2}].
	\end{equation}

	\subsubsection{Second-order version}
	
	The second-order approximation requires two iterations of (\ref{Duhamel}), that is 
	\begin{align}\label{variation2}
		u(h)=& e^{hA}u_0+ \int_0^he^{(h-t_1)A}B(t_1)e^{t_1A}u_0\,\text{d}t_1 \\&+ \int_0^{h}\int_{0}^{t_1}e^{(h-t_1)A}B(t_1)e^{(t_1-t_2)A}B(t_2)e^{t_2A}u_0\,\text{d}t_2\text{d}t_1 + R_{V_2},\non
	\end{align}
	where the $\mathcal{O}(h^3)$ remainder is given by 
	\begin{align}\label{RV2}
		R_{V_2}=&\int_0^{h}\int_{0}^{t_1}\int_{0}^{t_2}e^{(h-t_1)A}B(t_1)e^{(t_1-t_2)A}B(t_2)e^{(t_2-t_3)A}B(t_3)u(t_3)\,\text{d}t_3\text{d}t_2\text{d}t_1.
	\end{align}
	In order to derive the family of integrators (\ref{family}), we employ again the quadrature (\ref{quadrature1}) for the single integral in (\ref{variation2}). In turn, the double integral must be approximated according to the following rule
	\begin{equation}\label{quadrature2}
		I_2:=\  \frac{h^2}{8}f_2(h(1-\ta),h(1-\ta)) + \frac{h^2}{4}f_2(h(1-\ta),h\ta) + \frac{h^2}{8}f_2(h \ta,h \ta),
	\end{equation}	
	where 
	\begin{equation}
		\label{f2}
		f_2(t_1,t_2)= e^{(h-t_1)A}B(t_1)e^{(t_1-t_2)A}B(t_2)e^{t_2A}u_0. 
	\end{equation}

	\begin{rem}
		Note that for $\tau=\tfrac{1}{2}$, the  quadrature $I_2$ indicates that the inner integral is approximated with the right-point rule, while the outer one is  via the midpoint one. For $\tau\in[0,\tfrac{1}{2})$ the quadrature's $I_2$ nodes are distributed as shown in Fig. \ref{fig:RD}.

		\begin{figure}[h]
			\centering
			\includegraphics[]{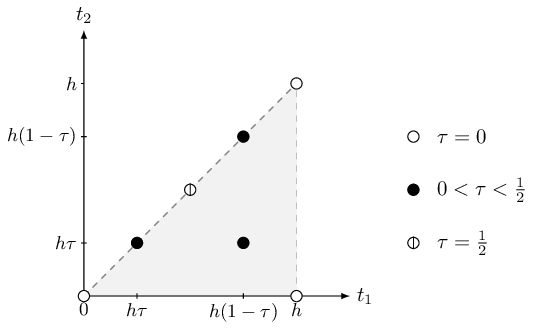}
			\caption{Distribution of the nodes in quadrature $I_2$ over the 2-d simplex (triangle), see (\ref{quadrature2}). At $\ta=0$ they lie on the vertices of the triangle. Meanwhile, at $\tau=1/2$ they all coincide at $(h/2,h/2)$. For both cases, they define a boundary-type quadrature and result in splitting featuring  three exponentials. For $0<\tau<1/2$, we arrive at methods composed of five exponentials, and the  nodes lie  as  displayed in the figure.}
			\label{fig:RD}
		\end{figure}
	\end{rem}
	The error term of the quadrature $I_2$
	\begin{equation}\label{RI2}
		R_{I_2}:=\int_0^{h}\int_0^{t_1}f_2(t_1,t_2)\,\text{d}t_2\text{d}t_1\ -\  I_2,
	\end{equation}	
	will be presented in Section \ref{error}. As previously, we observe that 
	\begin{align}\label{still_missing}
		e^{hA}u_0+I_1+I_2 = &\, {\rm e}^{h\tau A} \,  {\rm Id}\,  e^{h(1-2\tau) A} \,  {\rm Id} \,  {\rm e}^{h\tau A}u_0\non \\
		+&\frac{h}{2}\,e^{ h\tau A}\,  B(h(1-\tau))\,  e^{h(1-2\tau) A}\,  {\rm Id}\,  e^{h\tau A}u_0\non \\ 
		+&\frac{h}{2}\,e^{h\tau A}\,  {\rm Id}\,  e^{h(1-2\tau) A}\,  B(h\tau)\,  e^{h\tau A}u_0 \non \\ 
		+&\frac{h^2}{8}\,e^{ h\tau A}\,  B^2(h(1-\tau))\,  e^{h(1-2\tau) A}\,  {\rm Id}\,  e^{h\tau A}u_0\non \\
		+&\frac{h^2}{8}\,e^{h\tau A}\,  {\rm Id}\,  e^{h(1-2\tau) A}\,  B^2(h\tau)\,  e^{h\tau A}u_0 \non \\
		+&\frac{h^2}{4}e^{ h\tau A}\,  B(h(1-\tau))\,  e^{h(1-2\tau)A}\,  B(h\tau)\,  e^{h\tau A}u_0.
	\end{align}
	With three additional summands,
	\begin{align}
		R_{E_2}:=&\, \frac{h^3}{8}\,e^{h\ta A}\,  (B(h(1-\ta)) \,  e^{h(1-2\ta)A} \,  B^2(h\ta)\,  e^{h\ta A}u_0\non \\
		+\, &\, \frac{h^3}{8}\,e^{h\ta A}\,  B^2(h(1-\ta))\,  e^{h(1-2\ta)A}\,  B(h\ta)\,  e^{h\ta A}u_0\non \\
		+\, &\,  \frac{h^4}{16}e^{h\ta A}\,  B^2(h(1-\ta))\,  e^{h(1-2\ta)A}\,  B^2(h\ta)\,  e^{h\ta A}u_0,
		\label{RE2}
	\end{align} 
	the approximation (\ref{still_missing}) leads further to
	\begin{align*}
		e^{hA}u_0+I_1+I_2+R_{E_2} =e^{ h\tau A}&\left[{\rm Id} +  \frac{h}{2} B(h(1-\tau))+  \frac{h^2}{8} B^2(h(1-\tau)) \right]e^{h(1-2\tau)A}\times\\
		&\left[{\rm Id}\, +\,  \frac{h}{2}\,B(h\tau) + \frac{h^2}{8} B^2(h(1-\tau)) \right]e^{h\tau A}u_0.
	\end{align*}
	By similar arguments as in the previous subsection we can replace the terms inside parenthesis by semigroups. This replacement introduces another source of error, namely
	\begin{align}
		R_{S_2}:=&\frac{h^6}{48^2}e^{ h\tau A}  B^3(h(1-\tau))e^{\frac{h}{2}\xi_1B(h(1-\tau)}e^{h(1-2\tau)A} B^3(h\tau) e^{\frac{h}{2}\xi_2B(h\tau)}e^{h\tau A}u_0 \non \\
		-&\frac{h^3}{48}e^{ h\tau A}e^{\frac{h}{2} B(h(1-\ta))}e^{h(1-2\tau)A} B^3(h\tau)e^{\frac{h}{2}\xi_2B(h\tau)} e^{h\tau A}u_0 \non \\
		-&\frac{h^3}{48}e^{ h\tau A} B^3(h(1-\tau))e^{\frac{h}{2}\xi_1B(h(\tau-1))} e^{h(1-2\tau)A} e^{\frac{h}{2}B(h\ta) } e^{h
			\tau A}u_0,
		\label{RS2}
	\end{align}
	for certain $\xi_1\in[0,1]$ and $\xi_2\in[0,1]$. Finally, we conclude that
	\begin{equation}\label{Error_2}
		u(h)\ = \ \mathcal{F}(h,\ta)u_0+R_{V_2}+R_{I_1}+R_{I_2}-R_{E_2}+R_{S_2},\qquad \tau\in[0,\tfrac{1}{2}].
	\end{equation}
	
	\subsection{Error terms and their bounds}\label{error}
	In this section we establish the error bound of the error terms involved in the derivation of the family.
	Error terms $R_{E_1}$, $R_{S_1}$, $R_{E_2}$ and $R_{S_2}$  given by formulas (\ref{RE1}), (\ref{RS1}), (\ref{RE2}) and (\ref{RS2}) lead to straightforward estimates,
	\begin{align}
		\norm{R_{E_1}}\ \leq\ h^2 C_{E_1}\qquad {\rm and}\qquad \norm{R_{S_1}}\ \leq\ h^2 C_{S_1},\ \\
		\norm{R_{E_2}}\ \leq\ h^3 C_{E_2}\qquad {\rm and}\qquad \norm{R_{S_2}}\ \leq\ h^3 C_{S_2},\ 
	\end{align}
	where  error constants depend on $C_h$, $\|B(\cdot)\|_Y$ and $\|u_0\|$.
	On the other hand, bounds for  $R_{V_1}$, $R_{V_2}$, $R_{I_1}$, $R_{I_2}$, together with their expressions, require more attention. For that reason, we devote a subsection to each one.

	\subsubsection{Error  terms $R_{V_1}$ and $R_{V_2}$}
	
	Formulas (\ref{RV1}) and (\ref{RV2}) result directly in estimates
	
	\begin{equation*}
		\|R_{V_1}\|\leq\frac{h^2}{2!}\, C^2_h \|B\|^2_Y\,\|u\|_Y\quad {\rm and} \quad
		\|R_{V_2}\|\leq\frac{h^3}{3!}\, C^3_h \|B\|^3_Y\,\|u\|_Y,
	\end{equation*}
	which depend on the norm of the solution $\|u\|$ on the interval $[0,h]$. It is convenient, however, to express the estimates of $R_{V_1}$ and $R_{V_2}$ in terms of $\|u_0\|$, which is known. This can be done once $R_{V_1}$ and $R_{V_2}$ are presented as an infinite series of  nested integrals.
	To do so,  let us define the following operator
	$$
	K[u](h):=\int_0^h\,e^{(h-s)A}B(s)u(s)\,\text{d}s,
	$$
	and consider its $d$-multiple nesting  
	$$
	K^d[{\rm e}^{\cdot A}u_0](h)=\int_{0}^{h}\ldots\int_{0}^{t_{d-1}}{\rm e}^{(h-t_1)A}B(t_1)\ldots {\rm e}^{(t_{d-1}-t_d)A}B(t_{d})\,{\rm e}^{t_dA}u_0 \,\text{d}t_d\ldots\text{d}t_1,
	$$
	where $K^0[{\rm e}^{\cdot A}u_0]:={\rm e}^{h A}u_0$. Now, we can observe that iterations of (\ref{Duhamel}) lead to the Neumann series,
	$$
	u^{[n]}(h)\ =\ {\rm e}^{hA}u_0+K[u^{[n-1]}](h)\ =\ \sum_{d=0}^nK^d[{\rm e}^{\cdot A}u_0](h),\qquad n\geq1.
	$$
	This series converges to the solution of (\ref{main}), that is $$\displaystyle{u(h):=\lim_{n\rightarrow\infty}u^{[n]}(h)=\sum_{d=0}^\infty\,K^d[{\rm e}^{\cdot A}u_0](h)}.$$ Indeed, one can easily check that
	\begin{gather*}
		\|K^d[{\rm e}^{\cdot A}u_0]\|_Y\leq\\ C_{h}^{d+1}\,\|B\|_Y^d\,\|u_0\|\,\int_{0}^{h}\int_{0}^{t_1}\int^{t_2}_0\ldots\int_{0}^{t_{d-1}} \,\text{d}t_d\ldots\text{d}t_2\text{d}t_1 \leq \frac{(h\,C_{h}\,\|B\|_Y)^d}{d!}C_{h}\,\|u_0\|,
	\end{gather*}
	and  that
	$$
	\sum_{d=0}^\infty\,\|K^d[{\rm e}^{\cdot A}u_0]\|_Y\leq {\rm e}^{(h\, C_{h}\,\|B\|_Y)}\,C_{h}\,\|u_0\|.
	$$
	As a result, we conclude that $\{u^{[n]}\}_{n=0}^\infty$ forms a Cauchy sequence. To show that its limit solves (\ref{main}), we recall that $K^0[{\rm e}^{\cdot A} u_0]= {\rm Id}\, {\rm e}^{h A}u_0$ and observe 
	$$
	({\rm Id}-K)[u](h)=({\rm Id}-K)\left[ \sum_{d=0}^\infty\,K^d[{\rm e}^{\cdot A}u_0](h)\right]=(K^0-\lim_{n\rightarrow\infty} K^n)[{\rm e}^{\cdot A}u_0](h)={\rm e}^{h A}u_0.
	$$
	Now, it is obvious that $u(h)-u^{[n]}(h)=K^{n+1}\sum_{d=0}^\infty K^d[{\rm e}^{\cdot A}u_0](h)$ and that
	$$
	\|u(h)-u^{[n]}(h)\|\leq C_{h}{\rm e}^{(h\,C_{h}\,\|B\|_Y)}\,\frac{(h\,C_{h}\,\|B\|_Y)^{n+1}}{(n+1)!}\,\|u_0\|.
	$$
	Thus,  for the case of our interest $n=1,2$, we have
	\begin{align*}    
		\|R_{V_1}\|=&\|u(h)-u^{[1]}(h)\|\leq {\rm e}^{h\,C_h\,\|B\|_Y}\,\frac{(h\, C_h\,\|B\|_Y)^{2}}{2!}C_h\,\|u_0\|,\\
		\|R_{V_2}\|=&\|u(h)-u^{[2]}(h)\|\leq {\rm e}^{h\,C_h\,\|B\|_Y}\,\frac{(h\, C_h\,\|B\|_Y)^{3}}{3!}C_h\,\|u_0\|,
	\end{align*}
	respectively.
	\subsubsection{Quadrature error terms: $R_{I_1}$ and $R_{I_2}$}\label{QuadratureErrors}
	The main advantage of deriving the family of integrators (\ref{family}) is that neither the boundedness of $A$ nor of the commutators involving $A$  are required. However, error terms arising from quadratures involve commutators of $A$, $B(\cdot)$ and $B'(\cdot)$ as it will be presented in this subsection. Thus, to give sense to the bounds of the derived error terms, we make the following assumption
	\begin{assumption}\label{commutators_assum}
		Let us assume that 
		\begin{align*}
			(a)&\quad  \left\|[B(s),A]e^{rA}u_0\right\|,\quad \left\|[B'(s),A]e^{rA}u_0\right\|,\\
			(b)&\quad \left\|[[B(s),A],A]e^{sA}u_0\right\|,
		\end{align*}
		are well defined and bounded for $r,s\in[0,h]$.
	\end{assumption}
	We are now ready  to  present the derivation of the error terms arising from quadratures $I_1$ and $I_2$. The explicit form of the error terms will be relevant to find an appropriate value of $\tau$ that may  lead to the smallest error under additional assumptions on the  commutators involved in \ref{commutators_assum}. The approach we use does not employ Taylor series, which may lead to additional considerations when unbounded operators are involved. Instead, we use standard one-dimensional integration by parts. The key observation is that the boundary terms  lead to the quadrature rules as long as the proper choice of constants of integration  is made. In turn, the error of the quadrature emerges from the remaining integrals. We present a detailed derivation of $R_{I_1}$, and provide the final formula for $R_{I_2}$, which was obtained using the same approach. We skip its derivation due to compactness of presentation.
	
	Consider  the integral of $f_1(t_1)\in C^2([0,h])$
	\begin{equation*}
		\int_0^{h} f_1(t_1)\,\text{d}t_1\  =\ h \int_0^{1} f_1(h s)\,\text{d}s,\qquad s=\frac{t_1}{h},
	\end{equation*}
	and split the integral of the right-hand side as follows
	\begin{equation*}
		\int_0^{1} f_1(h s)\,\text{d}s\ = \ \int_0^{\ta} f_1(h s)\,\text{d}s\ +\ \int_\tau^{1-\tau} f_1(h s)\,\text{d}s\ +\ \int_{1-\ta}^{1} f_1(h s)\,\text{d}s.
	\end{equation*}
	Integrating by parts each integral,  we arrive at
	\begin{align*}
		\int_0^{1} f_1(h s)\,\text{d}s\ =\ &-a_1f_1(0)+ (a_1-b_1)f_1(h\tau) + (b_1-c_1)f_1(h(1-\ta)) + (1+c_1)f_1(h)\non\\
		&-\int_0^{\ta}(s+a_1)\left(\frac{\pa f_1(h s)}{\pa s}\right)\text{d}s -\int_\ta^{1-\ta}(s+b_1)\left(\frac{\pa f_1(h s)}{\pa s}\right)\text{d}s \non\\
		&-\int_{1-\ta}^{1}(s+c_1)\left(\frac{\pa f_1(h s)}{\pa s}\right)\text{d}s,
	\end{align*}
	where $a_1$, $b_1$, and $c_1$ are arbitrary constants of integration. 
	To reconstruct the quadrature $I_{1}$ from the  boundary terms, 
	we choose $a_1=0$, $b_1=-1/2$, and $c_1=-1$, cf. (\ref{quadrature1}). Therefore,
	\begin{equation}
		\label{1stform}
		\int_0^{h} f_1(t_1)\,\text{d}t_1 =\  \frac{h}{2}f_1(h\tau)\ +\ \frac{h}{2}f_1(h(1-\ta))\ +\ R^{(1)}_{I_1},
	\end{equation}
	where
	\begin{align}\label{RI11}
		R^{(1)}_{I_1} =&-h\int_0^{\ta}s\left(\frac{\pa f_1(h s)}{\pa s}\right)\text{d}s -h\int_\ta^{1-\ta}\left(s-\frac{1}{2}\right)\left(\frac{\pa f_1(h s)}{\pa s}\right)\text{d}s \\
		&-h\int_{1-\ta}^{1}(s-1)\left(\frac{\pa f_1(h s)}{\pa s}\right)\text{d}s \nonumber
	\end{align}
	is the error of the quadrature which guarantees  accuracy of order $\mathcal{O}(h^2)$ at least\footnote{Note that $\pa f_1(h s)/\pa s=hf_1'(hs)$.}. In order to obtain $\mathcal{O}(h^3)$ accuracy,  we integrate by parts $R^{(1)}_{I_1}$ in (\ref{1stform}), obtaining
	\begin{align*}
		R^{(1)}_{I_1} =\ &\ 
		h^2\,a_2f'(0)\ +\  h^2\left(-\frac{\ta}{2}\ -\ a_2+b_2\right)f'(h\ta )\non\\
		&+\  h^2\left(\frac{\ta-1}{2}-b_2+c_2\right)f'(h(1-\ta))\ +\  h^2\left(\frac{1}{2}-c_2\right)f'(h)
		\non \\
		&-h\int_0^{\ta}\left(\frac{s^2}{2}+a_2\right)\left(\frac{\pa^2 f_1(h s)}{\pa s^2}\right)\text{d}s - h\int_\ta^{1-\ta}\left(\frac{s^2}{2}-\frac{s}{2}+b_2\right)\left(\frac{\pa^2 f_1(h s)}{\pa s^2}\right)\text{d}s\non\\
		&-h\int_{1-\ta}^{1}\left(\frac{s^2}{2}-s+c_2\right)\left(\frac{\pa^2 f_1(h s)}{\pa s^2}\right)\text{d}s,
	\end{align*}
	where $a_2$, $b_2$, and $c_2$ are new constants of integration. Choosing  $a_2=0$, $b_2=\ta/2$, and $c_2=1/2$,  boundary terms that involve derivatives are removed  and
	\begin{equation}
		\label{2ndform}
		\int_0^{h} f_1(t_1)\,\text{d}t_1 =\  \frac{h}{2}f_1(h\tau)\ +\ \frac{h}{2}f_1(h(1-\ta))\ +\ R^{(2)}_{I_1},
	\end{equation}
	where
	\begin{align}\label{RI12}
		R^{(2)}_{I_1}\  =\ & \frac{h}{2}  \int_0^{\tau }s^2 \left(\frac{\partial ^2f_1(h s)}{\partial s^2}\right)\text{d}s \ +\ \frac{h}{2}  \int_{\tau }^{1-\tau } \left(s(s-1)+\tau \right) \left(\frac{\partial ^2f_1(h s)}{\partial s^2}\right) \text{d}s\ \\
		&+\frac{h}{2} \int_{1-\tau }^1 (s-1)^2 \left(\frac{\partial ^2f_1(h s)}{\partial s^2}\right)\text{d}s. \non
	\end{align}
	The error term in quadrature (\ref{quadrature2}) can be derived in a similar way, i.e. based on integration by parts, resulting in
	\begin{gather}
		\label{RI2}
		R_{I_2}\ =\\
		-
		\frac{h^2}{2}\Int _{1-\tau }^1\!\Int _0^1\!\!\left(s_1^2-1\right)\! \left(\tfrac{\partial f_2(h s_1,h s_1 s_2)}{\partial s_1}\right)\!\text{d}s_2\text{d}s_1
		-
		\frac{h^2}{2}\Int _{\tau }^{1-\tau }\!\Int _0^1\!\!\left(s_1^2-\frac{1}{4}\right)\! \left(\tfrac{\partial f_2(h s_1,h s_1 s_2)}{\partial s_1}\right)\!\text{d}s_2\text{d}s_1\non\\ 
		-
		\frac{h^2}{2} \Int _0^{\tau }\Int _0^1s_1^2 \left(\tfrac{\partial f_2(h s_1,h s_1 s_2)}{\partial s_1}\right)\!\text{d}s_2\text{d}s_1
		-
		\frac{3h^2}{8}  \Int_0^{\frac{\tau }{1-\tau }} s_2 \left(\tfrac{\partial f_2(h (1-\tau ),h  (1-\tau )s_2)}{\partial s_2}\right)\!\text{d}s_2 \non \\ \nonumber
		- 
		\frac{3h^2}{8} \Int_{\frac{\tau }{1-\tau }}^1 \left(s_2-\frac{2}{3}\right) \left(\tfrac{\partial f_2(h (1-\tau ),h (1-\tau )s_2)}{\partial s_2}\right)\text{d}s_2 
		-  \frac{h^2}{8} \Int_0^1 s_2 \left(\tfrac{\partial f_2(h \tau ,h  \tau s_2 )}{\partial s_2}\right)\text{d}s_2.
	\end{gather}

	The  building blocks of $R^{(1)}_{I_1}$, $R^{(2)}_{I_1}$  are
	\begin{align}\nonumber
		\frac{\partial f_1(h s)}{\partial s}\ &=\  he^{h(1-s)A}\Big([B(hs),A]+B'(hs)\Big)e^{hsA}u_0\\ \label{nested_commutator}
		\frac{\partial ^2f_1(h s)}{\partial s^2}\ &=\ h^2e^{h(1-s)A}\Big(\left[B(hs),A],A\right]+2[B'(hs),A]+B''(hs)\Big)e^{hsA}u_0.
	\end{align}
	In turn, the building blocks for $R_{I_2}$ are 
	\begin{gather*}
		\frac{1}{h}\,\frac{\partial f_2(h s_1,h s_1 s_2)}{\partial s_1} =	\\ s_2 e^{h(1-s_1)A}B(hs_1)e^{hs_1(1-s_2)A}\Big([B(hs_1s_2),A]+ B'(hs_1s_2)\Big)e^{hs_1s_2A}u_0
	\end{gather*}
	\begin{gather*}
		\frac{1}{h(1-\tau)}\,\frac{\partial f_2(h (1-\tau )h (1-\tau )s_2)}{\partial s_2} = \\e^{h\ta A}B(h(1-\ta))\,e^{h(1-\ta)(1-s_2)A} \Big([B(h(1-\ta)s_2),A] + B'(h(1-\ta) s_2)\Big)e^{h\ta s_2A}u_0 
	\end{gather*}
	\begin{gather*}
		\frac{1}{h\tau}\,\frac{\partial f_2(h \tau ,h  \tau s_2 )}{\partial s_2} =\\ e^{h(1-\ta)A}B(h\ta)\,e^{h\ta(1-s_2)A}\Big([B(h\ta s_2),A] + B'(h\ta s_2)\Big)e^{h\ta s_2A}u_0.
	\end{gather*}
	Based on these formulas, we deduce the following bounds
	\begin{align}\non
		\|R^{(1)}_{I_1}\|\leq& h^2 C_a \max_{s\in[0,h]}\Big\{ \big\|[B(s),A]e^{sA}u_0\big\|,\ \|u_0\|\Big\}\\ \non
		\|R^{(2)}_{I_1}\|\leq& 
		h^3 C_b \max_{s\in[0,h]}\Big\{\big\|[B'(s),A]e^{sA}u_0\big\|,\ \big\| [[B(s),A],A]e^{sA}u_0\big\|,\ \|u_0\|\Big\}\\ \non
		\|R_{I_2}\|\leq& h^3 C_c \max_{r,s\in[0,h]}\Big\{\big\|[B(s),A]e^{rA}u_0\big\|,\ \| u_0\|\Big\}.
	\end{align}
	where the error constants have the following dependencies $C_a=C_a\big(\ta,C_h,\|B'(\cdot)\|_Y\big)$, $C_b=C_b\big(\ta,C_h,\|B''(\cdot)\|_Y\big)$ and $C_c=C_c\big(\ta,C_h,\|B(\cdot)\|_Y,\|B'(\cdot)\|_Y\big)$.

	\subsection{Accuracy of the family of integrators}\label{accuracy}

	A convenient feature of our approach is the simultaneous derivation and error analysis of the family of splittings $\mathcal{F}(h,\tau)$. Based on the above results, we  state the following Theorem.
	\begin{thm}
		Under Assumptions \ref{main_assum} and  \ref{commutators_assum}.a, the family of integrators (\ref{family})  satisfy 
		\begin{equation} 
			\|u(h)- \mathcal{F}(h, \tau)u_0\| \leq h^2\ C_1 \max_{s\in[0,h]}\Big\{\big\|[B(s),A]e^{sA}u_0\big\| , \|u_0\|\Big\},  \label{estimate1}
		\end{equation}
		meanwhile, under Assumptions\ \ref{main_assum} and \ref{commutators_assum}.b, it performs third (local) order of accuracy
		\begin{align} 
			\label{estimate2}
			\|u(h)- \mathcal{F}(h, \tau)u_0\| \leq h^3\ C_2  \max_{r,s\in[0,h]}\Big\{&\big\|[B(s),A]e^{rA}u_0\big\|, \big\|[B'(s),A]e^{sA}u_0\big\|, \\
			&\big\|[[B(s),A],A]e^{sA}u_0\big\|, \|u_0\|  \Big\},\non
		\end{align}
		where constant  $C_1$  depends on $\ta$,  $C_h$ and $\|B'\|_Y$, while $C_2$ additionally depends on $\|B\|_Y$ and $\|B''\|_Y$.
	\end{thm}

	The proof of the Theorem follows directly from the estimates of the error terms derived in the previous Sections. To conclude estimate (\ref{estimate1}), it is enough to observe that it comes from formula (\ref{Error_1}) and consists in terms: 
	$R_{V_1}$, $R_{I_1}$, $R_{E_1}$ and $R_{S_1}$; see (\ref{RV1}), (\ref{RI11}), (\ref{RE1}) and (\ref{RS1}), respectively. All terms but  $R_{I_1}$ are natural and easy to handle by most of problems. As mentioned in Subsection \ref{QuadratureErrors},  the quadrature error $R_{I_1}$ leads to an error constant with dependence on $\max_{s\in[0,h]}\|[B(s),A]e^{sA}u_0\|$. To establish estimate (\ref{estimate2}), we consider the relation (\ref{Error_2}) and appearing there error terms: $R_{V_2}$, $R_{I_1}$, $R_{I_2}$, $R_{E_2}$ and $R_{S_2}$; see (\ref{RV2}), (\ref{RI12}), (\ref{RE2}) and (\ref{RS2}), respectively. Again, the main error term $\max_{r,s\in[0,h]}\Big\{\big\|[B(s),A]e^{rA}u_0\big\|, \big\|[B'(s),A]e^{sA}u_0\big\|, \big\|[[B(s),A],A]e^{sA}u_0\big\|  \Big\}$\\ comes from the quadrature defects $R_{I_1}$ and $R_{I_2}$ defined with (\ref{RI12}) and (\ref{RI2}), respectively.

	\subsection{The role of parameter $\tau$}
	\label{role}
	The choice of $\tau\in[0,\frac{1}{2}]$  determines the number of exponents of the integrators. As mentioned before,  $\tau=0$ or $\tau=\frac{1}{2}$ yields  only three exponentials that need to be computed at each time step. All other $\tau\in(0,\frac{1}{2})$ make the computations more costly since the number of exponentials raises to five.  In this case, CPU times is $\sim1.7$ larger compared to the time required by $\tau=0,\frac{1}{2}$. Secondly, the value of $\tau$  may influence the accuracy of the approximation, but not its order. It is important to emphasize that there is no universal value of parameter $\tau$ that minimizes the error of the method. To explain this, let us observe that  error terms depending on $\tau$ are the following: $R_{E_1}$, $R_{S_1}$, $R_{E_2}$, $R_{S_2}$, $R_{I_1}$ and $R_{I_2}$. The aforementioned error terms involve: $B(s)$, $B'(s)$, $B''(s)$, ${\rm e}^{h B(s)}$, ${\rm e}^{ h  A}$, $[B(s),A]$, $[B'(s),A]$, $[[B(s),A],A]$.  Due  to the structure of the error terms, looking for the optimal $\tau$ is nearly impossible even if  $A$, $B$ and $u_0$ are known. 
	However, we can investigate the influence of $\tau$ on the accuracy of the method under additional assumptions. Let us consider the second-order convergence of $\mathcal{F}(h,\tau)$ for a general case of problem (\ref{main}) and assume that for each $s\in [0,T]$ the norm of the commutator $\|[[B(s),A],A]{\rm e}^{sA}u_0\|$ dominates  over the norms of all other error terms depending on $\tau$ (listed above).  This assumption is frequently fulfilled, especially when $A$ is a differential operator.  Considering a sufficiently small time step of integration, local error terms depending on $\tau$ are  slowly varying functions in each time-step and, without loss of generality, we may assume that they have a constant sign (in each time  interval). Thus, we may expect that  $\tau\approx0.21$ minimizes the error. To justify this value, let us observe that the dominating part $\|[[B(s),A],A]{\rm e}^{sA}u_0\|$ arises from $f_1''(hs)$, see (\ref{nested_commutator}), which in turn appears in $R_{I_1}^{(2)}$. The latter error is conveniently written in integral form.
	\begin{equation*}
		R_{I_1}^{(2)}\ =\ \int_0^1 K_{I_1}^{(2)}(s;\tau) f_1''(hs)\,\text{d}s,
	\end{equation*}
	with the $\tau$-dependent kernel
	\begin{equation}
		\label{kernel}
		K_{I_1}^{(2)}(s;\tau)\ =\ \frac{h^3}{2}\cdot \left\{
		\begin{aligned}
			s^2, &\qquad s \in[0,\tau]\\
			s(s-1)+\tau, &\qquad s \in[\tau,1-\tau]\\
			(s-1)^2, &\qquad s \in[1-\tau,1].
		\end{aligned}
		\right.
	\end{equation}
	By  demanding
	\begin{equation}
		\label{criterium}
		\int_0^1 K_{I_1}^{(2)}(s;\tau)\,\text{d}s\ =\ 0 ,
	\end{equation}
	we find the appropriate parameter  $\tau$ for which the error is reduced.
	By solving (\ref{criterium}), one may find that there is a single value of $\tau\in[0,1/2]$ with this property, leading to $\tau=\frac{1}{6} \left(3-\sqrt{3}\right)=0.211\,324...$. Given that the kernel $K_{I_1}^{(2)}(s;0.211\,324...)$ is multiplied by $\|f_1''(hs)\|\sim constant$  we may conclude that  the dominant error term  $\left|R_{I_1}^{(2)}\right|$ is minimized. Interestingly, this value of $\tau$ leads to the two-point Gauss-Legendre quadrature on the interval $s\in[0,1]$.
	

		
	
	{\color{black}
		\section{Family of Splittings $\mathcal{D}(h,\tau)$}\label{Family_D}
		
		In this Section we will focus on the derivation and error analysis of the {\em second-order version} of the  family of splittings $\mathcal{D}(h,\tau)$. We skip  details that were already discussed in Section \ref{Family_F}. Thus, we are only concerned with the global second-order performance of the family. Therefore, 
		additionally to Assumption (\ref{main_assum}), we also consider 
		\begin{assumption}\label{D_assum}
			We assume that commutator $[B(s),A]$ is bounded for each $s\in[0,h]$. 
		\end{assumption}
		This assumption is fulfilled when $A\sim\partial_x$ and $B(\cdot)$ is a multiplication operator, like in transport and Dirac equations, for which  $[B(s),A]$ is a multiplication operator.}
	
	{\color{black}
		\subsection{Derivation of the family of splittings} 
		The starting point is the twice-iterated Duhamel's formula (\ref{variation2}), where the single and double integrals are approximated by the quadratures 
	\begin{equation}
		\label{newQ}
		I_1:=\  hf_1(h\ta)+\frac{h^2(1-2\tau)}{2}f_1'(h\tau),\qquad I_2:=\  \frac{h^2}{2}f_2(h\ta,h\tau),
	\end{equation} 
	respectively, where $f_1$ and $f_2$ are defined like in the previous Section.
	In contrast to (\ref{quadrature1}), $I_1$  is now a Birkhoff quadrature due to the presence of $f_1'(h\tau)$. The underlying idea behind $\mathcal{D}(h,\tau)$ is  simple: instead of taking two  nodes to sample the integrand (\ref{variation1}) to reach $\mathcal{O}(h^3)$ in the accuracy of $I_1$, we can introduce a derivative evaluation in (\ref{newQ}). 
	Quadratures (\ref{newQ}) are the building blocks behind the family of splittings $\mathcal{D}(h,\tau)$. After the substitution of $f_1$ and $f_2$, see (\ref{quadrature1}) and (\ref{f2}), we find that 
	
	\begin{align*}
		e^{hA}u_0+I_1+I_2 =e^{ h(1-\tau) A}&\left[{\rm Id} + hB(h\tau)+ \frac{h^2(1-2\tau)}{2} \left( [B(h\tau),A]+B'(h\tau)\right)\right. \\
		&\left.+\frac{h^2}{2}B^2(h\tau)  \right]e^{h\tau A}u_0.
	\end{align*}
	Thanks to the Assumptions (\ref{main_assum}) and (\ref{D_assum}), the expression inside the square parenthesis can be approximated  by $e^{hB(h\tau)+\frac{h^2(1-2\tau)}{2}\mathcal{C}(h\tau)}$ with $\mathcal{C}(t)= [B(t),A]+B'(t)$. Following the same procedure as in the previous section, we end up with\begin{equation}
		\label{D-error}
		u(h)\ = \ \mathcal{D}(h,\ta)u_0+R_{V_2}+R_{I_1}+R_{I_2}+R_{S},\qquad \tau\in[0,1].
	\end{equation}
	\begin{align*}
		R_S=-e^{h(1-\tau) A}&\left(\frac{h^3(1-2\tau)}{4}\left(B(h\ta)\mathcal{C}(h\tau)+\mathcal{C}(h\tau)B(h\ta)\right)+\frac{h^4(1-2\tau)^2}{8}\mathcal{C}^2(h\tau)\right. \\
		&\left.+\frac{1}{6}\left[hB(h\tau)+\frac{h^2(1-2\tau)}{2}\mathcal{C}(h\tau)\right]^3\!\!e^{\xi (hB(h\tau)+\frac{h^2}{2}(1-2\tau)\mathcal{C}(h\tau))}\right)e^{h\tau A}u_0, 
	\end{align*}
	where $\xi\in[0,1].$ Correspondingly,
	\begin{equation*}
		\norm{R_S}\leq C_{S}h^3 ,   
	\end{equation*}
	where $C_S$ depends on $C_h$, $\norm{B(\cdot)}_Y$, $\norm{[B(\cdot),A]}_Y$ and $\norm{u_0}$. The definition of the other error terms in (\ref{D-error}), $R_{V_2}$, $R_{I_1}$, and $R_{I_2}$, can be found in   (\ref{RI1}), (\ref{RI2}) and (\ref{RV2}), respectively.

	\subsection{Error terms and their bounds}
	\label{family-D-error}
	Based on the same arguments used in Section \ref{QuadratureErrors}, we establish the error terms for each quadrature, $R_{I_1}$ and $R_{I_2}$, explicitly 
	\begin{equation}
		\label{RI2-D}
		R_{I_1}= \frac{h}{2}\int_0^\tau(s-1)^2\left(\frac{\pa^2f_1(hs)}{\pa s^2}\right)\text{d}s\ +\ \frac{h}{2}\int_\tau^1s^2\left(\frac{\pa^2f_1(hs)}{\pa s^2}\right)\text{d}s
	\end{equation}
	and
	\begin{align}
		R_{I_2}= -\frac{h^2}{2}&\left[\int_0^1\!\!s_2\left(\tfrac{\pa f_2(h\tau,h\tau s_2)}{\pa s_2}\right)\!\text{d}s_2+\int_0^\tau\!\!\int_0^1\!\!s_1^2\left(\tfrac{\pa f_2(h s_1,hs_1 s_2)}{\pa s_1}\right)\!\text{d}s_2\text{d}s_1\right.\\
		&+\left.\int_\tau^1\!\!\int_0^1\!\!(s_1^2-1)\left(\tfrac{\pa f_2(h s_1,hs_1 s_2)}{\pa s_1}\right)\!\text{d}s_2\text{d}s_1\right] .
	\end{align}
	Based on the above formulas, we deduce the  bounds for the quadrature errors
	\begin{align}\non
		\|R_{I_1}\|\leq& 
		h^3 C_b \max_{s\in[0,h]}\Big\{\big\|[B'(s),A]e^{sA}u_0\big\|,\ \big\| [[B(s),A],A]e^{sA}u_0\big\|,\ \|u_0\|\Big\}\\ \non
		\|R_{I_2}\|\leq& h^3 C_c \max_{r,s\in[0,h]}\Big\{\big\|[B(s),A]e^{rA}u_0\big\|,\ \| u_0\|\Big\}.
	\end{align}
	where the error constants have the following dependencies $C_b=C_b\big(\ta,C_h,\|B''(\cdot)\|_Y\big)$ and $C_c=C_c\big(\ta,C_h,\|B(\cdot)\|_Y,\|B'(\cdot)\|_Y\big)$.
	With these results at hand, we state the following result. 
	\begin{thm}
		Under Assumptions\ \ref{main_assum}, \ref{commutators_assum}.b and \ref{D_assum}, it performs third (local) order of accuracy
		\begin{align} 
			\|u(h)- \mathcal{D}(h, \tau)u_0\| \leq h^3\ C  \max_{r,s\in[0,h]}\Big\{&\big\|[B(s),A]e^{rA}u_0\big\|, \big\|[B'(s),A]e^{sA}u_0\big\|, \\
			&\big\|[[B(s),A],A]e^{sA}u_0\big\|, \|u_0\|  \Big\},\non
		\end{align}
		where constant  $C$  depends on $\ta$, $\|B\|_Y$ and $\|B''\|_Y$.
	\end{thm}
	Despite the appearance of the commutators and derivatives in the family $\mathcal{D}(h,\tau)$, we can notice the similarity in the results on the convergence of family $\mathcal{F}(h,\tau)$.
	\subsection{The role of parameter $\tau$}
	We now tackle the role of the parameter $\tau$ on the family $\mathcal{D}(h,\tau)$. The value ${\tau=1/2}$ is exceptional since it reduces the number of stages from three to two, leading to the midpoint Strang splitting, see (\ref{familyII}). There is no other value with this property. Under the same considerations discussed in Section \ref{role} about the dominance of $\|[[B(s),A],A]{\rm e}^{sA}u_0\|$, we may look for the optimal $\tau$ that minimizes its contribution. For this purpose,  we extract the kernel from (\ref{RI2-D}), namely 
	\begin{equation}
		\label{kernel}
		K_{I_1}^{(2)}(s;\tau)\ =\ \frac{h^3}{2}\cdot \left\{
		\begin{aligned}
			(s-1)^2, &\qquad s \in[0,\tau]\\
			s^2, &\qquad s \in[\tau,1].
		\end{aligned}
		\right.
	\end{equation}
	Note that $K_{I_1}^{(2)}(s;\tau)$ is a positive function, implying that its integral on $[0,1]$ is positive as well. As a consequence, there is no real $\tau$ with the property $\int_0^1K_{I_1}^{(2)}(s;\tau)\text{d}s=0$. However, $\tau=1/2$  minimizes $\int_0^1K_{I_1}^{(2)}(s;\tau)\text{d}s$.  Thus, the midpoint Strang splitting is potentially the best member of the family $\mathcal{D}(h,\tau)$. 
}

\section{Numerical Examples} \label{example}
\subsection{Schr\"odinger Equation}
To illustrate the performance of $\mathcal{F}(h,\tau)$, we consider the celebrated one-dimensional linear
Schr\"odinger equation in the following setting
\begin{equation}\label{Schrodinger}
	\left\{
	\begin{aligned}
		&i\frac{\partial u(x,t)}{\partial t} =\left(-\frac{1}{2}\frac{\partial^2 u(x,t)}{\partial x^2} + V(x,t)\right)u(x,t),\qquad t\in[0,1],\qquad x\in[-3,3],\\
		&u(x,0)\ =\ (x^2-3^2)\,e^{-20 (x+\frac{1}{2})^2},\\
		&u(\pm 3,t)=0,
	\end{aligned}
	\right.
\end{equation}
with the potential $V(x,t) =  -2\cos(10\,t)\,x^2+x^4$.
\begin{rem}
	Equation (\ref{Schrodinger}) can be considered as a special case of  (\ref{main}) once we assume that $u(\cdot,t)\in C^1([0,1],H^2[-3,3])$,  $A:H^2([-3,3])\rightarrow L^2([-3,3])$, $B(t)\in C^2([0,1],H^2[-3,3])$ and take
	\begin{equation}
		\frac{1}{i}A =-\frac{1}{2} \frac{\partial^2}{\partial x^2}\ ,\qquad \frac{1}{i}B(t)\ =\ V(x,t).
	\end{equation}
	In this setting operator $A$ is  unbounded, while $\|B(t)\|_{L^2[(-3,3)]}\leq 99$. Moreover, given the boundary conditions,  one can easily prove that $A$ and $B(t)$ are skew-Hermitian operators; therefore, $\|e^{hA}\|_{H^2[(-3,3)]}=\|e^{hB(t)}\|_{H^2[(-3,3)]}=1$ for each $t\in[0,1]$, and additionally  $\norm{u(\cdot,t)}_{L^2[(-3,3)]}\ =\ {\rm constant}$. In this setting, the performance of $\mathcal{F}(h,\tau)$ is expected to be of second-order (globally) obeying  (\ref{estimate2}). We are not discussing the first-order performance  $\mathcal{F}(h;\tau)$ since it was already addressed for the time-independent case in \cite{Lubich} for the Strang Splitting. This consideration can be straightforwardly extended for arbitrary $\tau$.
\end{rem}

\begin{figure}[h]
	\centering
	\includegraphics[width=0.7\textwidth]{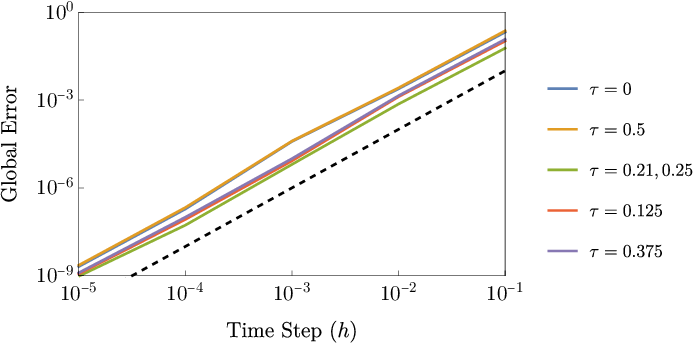}
	\caption {Global Error (in log-log scale) as a function of the time step $h$ for representative values of $\tau$. The black-dashed line represents the plot of $h^2$. The curves for $\tau=0.21,0.25$, see green curves, overlap each other.  }
	\label{globalerror}
\end{figure}

To carry out numerical experiments, we used the pseudo-spectral method described  in \cite{Baye}.  The spatial domain is discretized according to $N=250$ non-uniform mesh points $\{x_i\}_{i=1}^N$ distributed in $[-3,3]$ that serve as collocation points. They correspond to the zeros of  the  Legendre polynomial $L_{N}(\sigma x)$ with $\sigma=3$.
Then, via a time-stepping procedure, we evolved the initial condition from $t=0$ up to $t=1$. At the end of the evolution, at $t=1$, we compared the numerical solution with a reference one calculated by using a sufficiently small time step $h=10^{-6}$. Finally, we used the natural $L^2[(-3,3)]$ norm to compute the \textit{distance} between the reference solution and the approximate one. In Fig. \ref{globalerror}, we present in log-log  scale the plot of the global error  as a function of the time step $h$ for representative values of $\tau\in[0,\tfrac{1}{2}]$, paying special attention to the relevant value $\tau=0.21$.

Based on the numerical example and in accordance with the theoretical discussion about the influence of $\tau$ on $\mathcal{F}(h,\tau)$, we observed the following.
For all tested values of $\tau$, in particular $\tau=0,0.175,0.21,0.25,0.375,0.5$,  any member of the family $\mathcal{F}(h,\tau)$ performs second-order  error of accuracy (globally). It can be explicitly seen in Fig. \ref{globalerror}, where the plot of $h^2$ was included to compare with. In terms of accuracy, $\mathcal{F}(h,0.21)$ and $\mathcal{F}(h,0.25)$ lead to the smallest global error. We also considered the additional values  $\tau=0.01,0.2,0.22,0.49$ in calculations (not shown in Fig. \ref{globalerror} for simplicity of presentation), $\tau=0.20,0.22$ leads to accuracy similar to the one reached by $\tau=0.21,0.25$, respectively. \textcolor{black}{In principle,  $\mathcal{F}(h,0.21)$ gives optimal approximation according to the theoretical investigations in Subsection \ref{role}. To explain the optimal accuracy of splitting $\mathcal{F}(h,0.25)$ we remind, that it is equivalent to the composition of two midpoint Strang splittings with the time step $h/2$: $\mathcal{F}(h/2,0.5)$. This means that the error constant (of two Strang splittings with half the step) drops by the factor of 4, and is as computationally costly as other 5-stages splittings. Indeed, two Strang splittings per one time step $h$ result in five evaluations of exponents.} 

\textcolor{black}{We observe a continuous $\tau$-dependence on the accuracy of family $\mathcal{F}(h,\tau)$ and, in particular, a very mild dependence when $\tau\in[0.21,0,25]$. Comparing the \textit{optimized} five-stages members $\mathcal{F}(h,0.21)$ and $\mathcal{F}(h,0.25)$ with three-stages members $\mathcal{F}(h,0)$ and $\mathcal{F}(h,0.5)$, we observe that $\tau=0.21, 0.25$ provide more accurate results, specifically the global error is almost one order of magnitude smaller.  However, the choice of $\tau\in[0.21,0.25]$ is computationally costly. In fact, the computational time of five-stages integrators are about 1.7 times longer than the three-stages integrators when $\tau=0$ or $\tau=0.5$.
}

{\color{black}
	\subsection{Transport Equation}
	To illustrate the performance of the family $\mathcal{D}(h,\tau)$, we consider the one-dimensional (linear) transport equation with source/sink $f(x,t)$ given by
	\begin{equation}\label{Transport}
		\left\{
		\begin{aligned}
			&\frac{\partial u(x,t)}{\partial t}= -\frac{\partial u(x,t)}{\partial x}+f(x,t)u(x,t),\qquad t\in[0,T],\qquad x\in\mathrm{R},\\
			&u(x,0)\ =\ g(x),\\
			&u(\pm \infty,t)=0 .
		\end{aligned}
		\right.
	\end{equation}
	
	The exact solution of (\ref{Transport}) is known, namely
	\begin{equation}
		u(x,t)= \exp(\int _1^xf(s,t-x+s)\text{d}s-\int _1^{x-t}f(s,t-x+s)\,\text{d}s)g(x-t)\ .
	\end{equation}
	For concreteness, we choose $g(x)=e^{-2x^2}$ and $f(x,t)=-e^{-(2x-t)^2}$.
	\begin{rem}
		Equation (\ref{Transport}) can be considered as a special case of  (\ref{main}) once we assume that $u(\cdot,t)\in C^1([0,1],H^1[-\infty,\infty])$,  $A:H^{1}([-\infty,\infty])\rightarrow L^2([-\infty,\infty])$, $B(t)\in C^2([0,1],H^{1}[-\infty,\infty])$ and take
		\begin{equation}
			A =- \frac{\partial}{\partial x}\ ,\qquad B(t)\ =\ f(x,t). 
		\end{equation}
		In this setting, $A$ is unbounded,  $\norm{B(t)}_{L^2[-\infty,\infty]}\leq 1$ and ${\norm{B'(t)}_{L^2[-\infty,\infty]}\leq \sqrt{2/e}}$. Furthermore,  we have\\ ${\norm{[[B(t),A]]}_{L^2[-\infty,\infty]}\leq\sqrt{8/e}}$. Consequently, the performance of $\mathcal{D}(h,\tau)$ is expected to be of second order (globally) according to Section \ref{family-D-error}. 
	\end{rem}
	
	For this numerical example, we employ a centered finite-differences scheme of order $\mathcal{O}(\Delta x^4)$ with $\Delta x=0.002$ to ensure high \textit{spatial} accuracy in the numerical solution.    Then, via time-stepping on $\mathcal{D}(h,\tau)$, we evolve the initial condition from $t=0$ to $t=1$. In practice, we restricted the domain to $x\in[-3,4]$, which is sufficiently large to avoid reflections due to the artificial boundaries. At the end of the evolution we compared the numerical solution with the exact one using the $L^2[-\infty,\infty]$ norm.
	
	Figure  \ref{fig:global-D} shows the global error as a function of $h$ for representative values of $\tau\in[0,1]$.  
	\begin{figure}[h]
				\centering
		\includegraphics[width=0.7\textwidth]{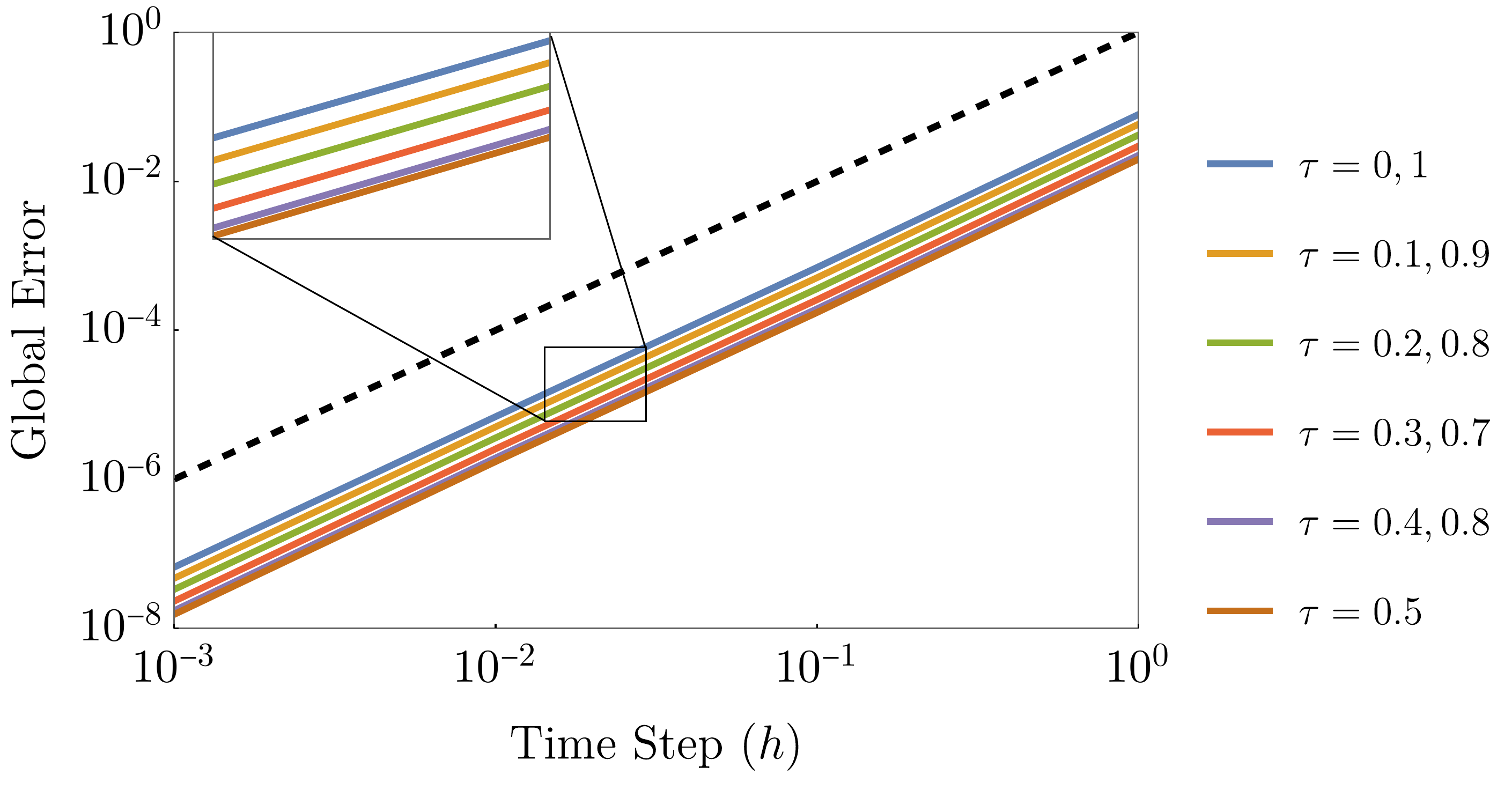}
		\caption{Global Error (in log-log scale) as a function of the time step $h$ for representative values of $\tau$. The black-dashed line represents the plot of $h^2$. The inset shows how $\tau=1/2$ leads to slightly better results than $\tau=0.4,0.8$.}
		\label{fig:global-D}
	\end{figure}
	We observe that $\tau=1/2$, leads to a smaller error compared with the other considered values. Thus, for this numerical example, $\mathcal{D}(h,1/2)$ can  regarded as the optimal member of the family. In addition, we observed that the global error for the members $\mathcal{D}(h,\tau)$ and $\mathcal{D}(h,1-\tau)$ is essentially the same: their plots overlap according to the line width we used, see Fig. \ref{fig:global-D}. This result is expected by inspecting  the structure of  (\ref{familyII}).
}

\section{Conclusions} \label{conclusions}
{\color{black}
	The quadratures used to approximate the integrals of an iterated Duhamel's formula are the building blocks to construct exponential splittings for linear differential equations of the form (\ref{main}). As illustration of this statement, we considered two families of second-order exponential splittings, namely 
	\begin{equation*}
		\mathcal{F}(h,\tau)= \  e^{h\ta A}e^{\frac{h}{2} B(h(1-\ta))}e^{h(1-2\ta) A}e^{\frac{h}{2}B(h\ta) }e^{h\ta A} ,\qquad \tau\in[0,\tfrac{1}{2}],
	\end{equation*}
	and
	\begin{equation*}  \mathcal{D}(h,\tau)= e^{h(1-\tau)A}e^{hB(h\tau)+\frac{h^2(1-2\tau)}{2}\mathcal{C}(h\tau)}e^{h\tau A},\quad\mathcal{C}(t)=[B(t),A]+B'(t),\quad\tau\in[0,1].
	\end{equation*}
	Despite the different structure of both families, their construction and error analysis can be done simultaneously in a constructive way using  the variation of constants formula and, in the most general case, Birkhoff quadratures. Our analysis includes the widely used midpoint Strang splitting, $\mathcal{F}(h,1/2)$ or/and $\mathcal{D}(h,1/2)$,  whose convergence was never proved in the case of unbounded operator $A$ and bounded, time-dependent operator $B(t)$. Although the idea of this extension was hinted in \cite{LubichBook},  the details were missing in the literature.
	
	Furthermore, we have illustrated the connection between exponential integrators and splitting methods in the following sense: splitting methods are specific types of exponential integrators, which requires higher regularity. Indeed, exponential integrators are based on (possibly iterated) variation of constants formula and  creative ways of numerical integration of the (possibly nested) integrals. We have shown that splitting methods may be derived via iterated variation of constants formula and that the quadratures of the integrals need specific nodes and weights.
	
	We have discussed the influence of parameter $\tau$ on $\mathcal{F}(h,\tau)$ and $\mathcal{D}(h,\tau)$ coming to the conclusion, that (in the sense of global accuracy) there is no optimal value of parameter $\tau$ suitable for all examples, but under additional assumptions the optimal value can be obtained by investigating the kernel associated to the one-dimensional quadrature error.

	One of the main advantages of the presented derivation and error analysis is that it can be straightforwardly extended to tackle higher-order splitting methods. This naturally requires more iterations of the variation of constants  formula and, more importantly, finding appropriate multivariate Birkhoff quadratures  that allow the reconstruction of an exponential splitting. It is worth noticing that higher-order {\em compact} splittings, also known as {\em splittings with modified gradient potential} methods, like those presented in \cite{Chin,BAYE2003337,PhysRevE.70.056703} are inside the scope of the present methodology. In these methods not only computation of ${\rm e}^{hA}$ and ${\rm e}^{hB(s)}$ are required, but also terms like ${\rm e}^{c_1hB(s)+c_2\,h^3[B(s),[A,B(s)]]}$. To generate the latter semigroup, we need to use quadratures over multidimensional simplices featuring  derivatives of the integrand. The extension of the present work to higher-order methods will be done in detail in a forthcoming communication.
	
}
\section*{Acknowledgments}
This work has been supported by The National Center
for Science (NCN), based on Grant No. 2019/34/E/ST1/00390.
Computations were carried out using the computers of Center of Informatics Tricity Academic Supercomputer \& Network (CI TASK).

\bibliographystyle{siamplain}
\bibliography{references}
\end{document}